\documentclass[12pt,twoside,reqno]{amsart}
\usepackage[colorlinks=true,citecolor=blue]{hyperref}
\usepackage{mathptmx, amsmath, amssymb, amsfonts, amsthm, mathptmx, enumerate, color,mathrsfs}
\usepackage{graphicx}

\usepackage{multirow}
\usepackage{epstopdf}
\usepackage{multicol}
\usepackage{algorithm}
\usepackage{algorithmic}
\usepackage{epstopdf}

\usepackage{nicefrac}       
\usepackage{cleveref}

\newcommand{\RR}{\mathbb{R}}

\newcommand{\cX}{\mathscr{X}}
\newcommand{\cY}{\mathscr{Y}}
\newcommand{\R}{\RR} 

\DeclareMathOperator*{\Id}{Id}
\newcommand{\cH}{{\mathcal H}}
\newcommand{\cL}{{\mathcal L}}
\newcommand{\cR}{{\mathcal R}}

\newcommand{\cU}{{\mathcal U}}

\newcommand{\cN}{{\mathcal N}}
\newcommand{\cZ}{{\mathcal Z}}

\newcommand{\T}{\mathscr{T}}
\newcommand{\tg}{\bar \mu}
\newcommand{\tgY}{\bar \nu}
\newcommand{\qEM}{\theta}
\newcommand{\Q}{\Theta}

\providecommand{\norh}[1]{\lVert{#1}\rVert}
\providecommand{\noridx}[2]{\lVert{#1}\rVert_{#2}}

\DeclareMathOperator*{\argmax}{argmax}
\DeclareMathOperator*{\argmin}{argmin}

\DeclareMathOperator*{\tr}{tr}
\DeclareMathOperator*{\KL}{KL}

\newtheorem{theorem}{Theorem}[section]

\newtheorem{proposition}[theorem]{Proposition}

\theoremstyle{definition}

\newtheorem{remark}[theorem]{Remark}
\numberwithin{equation}{section}

\begin{document}
\setcounter{page}{1}

\vspace*{2.0cm}
\title[Alternating minimization for estimation and identification]
{Alternating minimization for simultaneous estimation of a latent variable and identification of a linear continuous-time dynamic system}
\author[P-C. Aubin-Frankowski, A. Bensoussan, S. J. Qin]{Pierre-Cyril Aubin-Frankowski$^{1,*}$, Alain Bensoussan$^2$, S. Joe Qin$^3$}
\maketitle
\vspace*{-0.6cm}

\begin{center}
{\footnotesize

$^1$INRIA-Département d'Informatique de l'École Normale Supérieure,
PSL, Research University,Paris, France\\
$^2$International Center for Decision and Risk Analysis, Jindal School of Management, University of Texas at Dallas, School of Data Science, City University Hong Kong\\
$^3$Institute of Data Science, Lingnan University, Hong Kong

}\end{center}

\vskip 4mm {\footnotesize \noindent {\bf Abstract.}
We propose an optimization formulation for the simultaneous estimation of a latent variable and the identification of a linear continuous-time dynamic system, given a single input-output pair. We justify this approach based on Bayesian maximum a posteriori estimators. Our scheme takes the form of a convex alternating minimization, over the trajectories and the dynamic model respectively. We prove its convergence to a local minimum which verifies a two point-boundary problem for the (latent) state variable and a tensor product expression for the optimal dynamics.


 }

\renewcommand{\thefootnote}{}
\footnotetext{ $^*$Corresponding author.
\par
E-mail address: pierre-cyril.aubin@inria.fr (P-C. Aubin-Frankowski), axb046100@utdallas.edu (A. Bensoussan), sjoeqin@outlook.com (S. J. Qin).
\par
Received xx, x, xxxx; Accepted xx, x, xxxx.

\rightline {\tiny   \copyright  2023 Communications in Optimization Theory}}

\emph{In memory of Roland Glowinski,}

\section{Introduction}

The theory of latent variables in Data Science has been progressing
very fast in the recent years, with the objective of reducing the
dimension of the dataset. Since data is often associated with dynamic
systems, it is natural to consider in this context the framework of identification and estimation
of dynamic systems. We refer to \cite{Qin2020}
for a survey of the main ideas in this direction and to \cite{Yu2022}
for details, in connection with the Kalman filter in discrete time. The general
idea is to consider the latent variable as described by a dynamic
system in state space representation. The difficulty is that we
need to identify the system while estimating it. The maximum likelihood approach is a natural way to proceed. Here, our algorithm reads as an alternating minimization of a quadratic objective that is nonconvex due to a bilinear term, but convex component-wise. Alternating minimization, a special case of block coordinate descent, is among the simplest algorithms one may think of. However one of the authors has shown in \cite{leger2023gradient} that alternating minimization encompasses many algorithms (such as gradient descent and its variations). It is also well-known in statistics through the Expectation-Maximization (EM) algorithm \cite{Neal1998}, to which our algorithm corresponds. We refer to \cite[Section 4.8]{leger2023gradient} for more interpretations of EM. Nevertheless, because of the nonconvexity of the objective, we cannot expect the type of global convergence shown in \cite{leger2023gradient}. 

 A discrete time version for a simpler model of the algorithm can be found in \cite{Bensoussan2021}. While most of the papers on this specific topic are indeed in discrete time, there is a huge swath of literature dealing with dynamic systems in continuous time, justifying tackling this setting as well. To simplify the theory and the algorithm, we consider that some aspects of the linear dynamic system are known, in particular the covariance matrices of the noises and the observation
matrix. 

\section{The model}
We consider an input-output problem $(v(t),y(t))_{t\in[0,T]}$ in state space representation, where we want from a single trajectory to reconstruct the state and its dynamic equation, in other words to perform both estimation and system identification at the same time. We assume that the linear dynamic systems with noise are described by two stochastic differential equations, \emph{with all the underlined quantities being known},
\begin{align}
	dx&=(Ax+B \underline v)dt+\underline Gdw, \; &&x(0)=\xi\sim \cN(\underline x_{0},\underline \Pi_{0}) , \label{eq:2-1}\\
	dy&=\underline Cx(t)dt+db(t),\; &&y(0)=0, \label{eq:2-2}
\end{align}
where we assume for simplicity that the dimensions of the operators are fixed, for instance chosen minimal through realization theory, see e.g. \cite{Bensoussan2021}.  More precisely we take $x(t)\in \R^{N}$, $v(t)\in \R^{d}$, $A\in\cL(\R^{N};\R^{N})$, $B\in\cL(\R^{d};\R^{N})$;
$v(t)$ being a given deterministic control. There is no optimal control
in this setup, but $v(.)$ is an input decided by the controller.
The process $w(t)$ is a Wiener process in $\R^{m}$, with correlation
matrix $Q$, and $G\in\cL(\R^{m};\R^{N})$. The random variable $\xi$ is Gaussian with mean $x_{0}\in \R^{N}$ and covariance matrix
$\Pi_0$. It is independent of the Wiener process $w(t)$. The matrices
$A$ and $B$ are not known, although they are in the vicinity of known
matrices $A_{0},B_{0}$ used as priors. The state of the system $x(t)$ is not observable.
We observe instead the process $y(t)\in \R^{p}$ where $C\in\cL(\R^{N};\R^{p})$ is known. The process $b(t)$
is a Wiener process in $\R^{p}$, with covariance matrix $R$, also
independent from $\xi$ and $w(.)$. For
simplicity, we assume in the sequel that the matrices $\Pi_0,Q,R$ are invertible (inverses would then be replaced with pseudo-inverses). 


The only information is $v(t),y(t),t\in[0,T]$. If $A,B$ were known,
the problem would reduce to estimating the evolution of the state $x(t)$.
This is the classical Kalman smoothing problem. There are many equivalent
ways to solve it. For instance, it is well known that the maximum likelihood is
equivalent to the following least square deterministic control problem 
\begin{equation}
	\dfrac{dx}{dt}=Ax+Bv+Gw, \quad x(0)=\xi, \label{eq:2-3}
\end{equation}
in which the control is the pair $(\xi,w(.))$ minimizing the payoff\footnote{In continuous time, $y(t)\in L^2(\T,\R^m)$ ``is reminiscent of the observation process, in fact rather the	derivative of the observation process (which, as we know, does not exist)'' \cite[p180]{Bensoussan2018}. Thus it is as if we observed this derivative to do the reconstruction, which is justified since only integrals of it appear in \eqref{eq:2-4}.} 
\begin{equation}
	\dfrac{1}{2}\Pi_0^{-1}(\xi-x_{0}).(\xi-x_{0})+\dfrac{1}{2}\int_{0}^{T}Q^{-1}w(t).w(t)dt+\dfrac{1}{2}\int_{0}^{T}R^{-1}(y(t)-Cx(t)).(y(t)-Cx(t))dt.\label{eq:2-4}
\end{equation}
When $A,B$ are not known, one way to proceed is to approach (\ref{eq:2-3})
as a constraint. This leads to the following formulation. Introduce
the argument $Z=(A,B,x(.).w(.))$ where $x(.)\in H^{1}(0,T;\R^{N})$,
$w(.)\in L^{2}(0,T;\R^{m})$, and define the norm 
\begin{equation}
	||Z||^{2}=\tr(AA^{*}+BB^{*})+|x(0)|^{2}+\int_{0}^{T}|\dfrac{dx}{dt}|^{2}dt+\int_{0}^{T}|w(t)|^{2}dt, \label{eq:2-5}
\end{equation}
thus $Z$ belongs to a Hilbert space, denoted by $\cZ$. Based on \eqref{eq:2-4}, we then introduce the following functional $J$ over $\cZ$ to perform the reconstruction. It will be justified in \Cref{sec:bayesian} through Bayesian arguments.
\begin{align}
	\min_{(A,B,x(.),w(.))}& J((A,B,x(.),w(.)):=\dfrac{\alpha}{2}\tr\left((A-\underline A_{0})(A-\underline A_{0})^{*}+(B-\underline B_{0})(B-\underline B_{0})^{*}\right) \label{eq:2-6}\\
	&+\dfrac{\beta}{2}\int_{0}^{T}\left|\dfrac{dx}{dt}-Ax(t)-B\underline v(t)-\underline Gw(t)\right|^{2}dt+\dfrac{1}{2}\underline \Pi_{0}^{-1}(x(0)-\underline x_{0}).(x(0)-\underline x_{0})\nonumber \\
	&+\dfrac{1}{2}\int_{0}^{T}\underline Q^{-1}w(t).w(t)dt+\dfrac{1}{2}\int_{0}^{T}\underline R^{-1}(\underline y(t)-\underline Cx(t)).(\underline y(t)-\underline Cx(t))dt  \nonumber,
\end{align}
for given parameters $\alpha,\beta>0$. To alleviate notations, \emph{from now on we do not underline the known quantities} and we refer the reader to this section. The first term in $\alpha$ is a regularizing term in $(A,B)$ using our prior. The second term in $\beta$ is a penalty term if the constraint (\ref{eq:2-3}) is not satisfied. It can be seen as a regularization of \eqref{eq:2-4} (see \Cref{rmk:inf_lim} below). The other terms penalize deviations of $(x(0),w(t),Cx(t))$ from their references $(x_0,0,y(t))$. So the problem amounts to minimizing the functional $J(Z)$ on the Hilbert space $\cZ$. Note that $(A,B)$ only appear in the two first terms, while $(x(.),w(.))$ appear in all terms but the two first. This suggests to do an alternating minimization of $J$ as in \Cref{sec:algo} below. However $J$ is non convex due to the term $Ax(t)$, so we cannot hope to reach a global minimum for every initialization. We will first justify the choice of $J$ in \Cref{sec:bayesian} and then prove the existence of a minimum in \Cref{sec:existence}, giving also the first-order optimality conditions that it satisfies. Note that a similar methodology can be replicated if we consider some other matrices to be unknown (e.g.\ $C$ or $G$). In the discrete time case, the problem of estimating simultaneously two matrices while minimizing an expression of their product has been considered in \cite{Qin2022}.

\section{Bayesian justification of the model}\label{sec:bayesian}

We follow here the presentation of \cite[Section 2]{Guth2022}  for the classical derivation of a least squares problem from a maximum a posteriori (MAP) estimator of system based on a model ($M$) and observations operator ($O$). To avoid technicalities, we do the justification for random variables over a finite set, thus not for the stochastic differential equations \eqref{eq:2-1}-\eqref{eq:2-2} with Brownian motions which we considered. However the ideas and results extend to infinite dimensions \cite{Dashti2013}.\footnote{The key reason for the methodological difference when moving to continuous time is reminded in \cite{Dashti2013} ``While in the finite-dimensional
	setting, the prior and posterior distribution of such statistical problems can typically
	be described by densities w.r.t. the Lebesgue measure, such a characterisation is no
	longer possible in the infinite dimensional spaces [...] no analogue of the Lebesgue measure exists in infinite dimensional spaces.'' However Gaussian measures can still serve as a replacement in our case \cite{DaPrato2006}. They correspond here to the Wiener processes we consider. Radon--Nikodym derivatives are then obtained through Girsanov's theorem \cite[Chapter 6.5]{Bensoussan2018}.} To obtain the objective function \eqref{eq:2-6}, we consider that \eqref{eq:2-2} is an equation of the form
\begin{equation*}
	Y=O(x(.),w(.))+\eta_{obs}
\end{equation*}
where $Y=(y_t)_{t\in \T}$, $\eta_{obs}\sim \cN(0,\cR)$, $\cR(s,t)=\delta_{s=t} R(t)$. Similarly we relax \eqref{eq:2-3} by introducing a model error
\begin{equation*}
	0=M(x(.),w(.),A,B)+\eta_{model}
\end{equation*}
where $\eta_{model}\sim \cN(0,\Id/\beta)$. We put a Gaussian prior $\mu_0$ on $(x(.),w(.),A,B)$ of the form $x(0)\sim\cN(x_{0},\Pi_{0})$, $w(.)\sim\cN(0,\mathcal Q)$ with $\mathcal Q(s,t)=\delta_{s=t} Q(t)$, $A\sim \cN(A_0,\Id/\alpha)$ and $B\sim \cN(B_0,\Id/\alpha)$. Thus, by Bayes' theorem, the posterior distribution is given by 
\begin{equation*}
	\mu^*(dZ)\propto \exp(-\frac{\beta}{2}\norh{M(x(.),w(.),A,B)}^2-\frac{1}{2}\noridx{y(.)-O(x(.),w(.))}{\cR})^2\mu_0(dZ).
\end{equation*}
where $\noridx{y(.)}{\cR}^2=\mathcal (\cR^{-1}y(.)).y(.)$. The MAP estimator is then given by $\argmax_{Z\in \cZ} \mu^*(dZ)$, and thus equivalently by minimizing the log-density, $\argmin_{Z\in \cZ} -\log \mu^*(dZ)$ which is precisely \eqref{eq:2-6}. More formally, in continuous time, to derive \eqref{eq:2-6} as the problem solved by the MAP, one can just apply \cite[Corollary 3.10]{Dashti2013} to identify $J$ as an Onsager-Machlup functional.

\begin{remark}[RKHS constraint in the limit case]\label{rmk:inf_lim} Interestingly \cite[Proposition 1]{Guth2022} recalls that for $\beta \rightarrow \infty$ (vanishing model noise case), we have that the accumulation points of the optimum $\hat Z_\beta$ all satisfy \eqref{eq:2-3}. For given $A,B$, \eqref{eq:2-3} says that $x(.)$ belongs to the affine vector space of functions $\cH=\{x(.) \, | \, \exists w(.), \, \nicefrac{dx}{dt}=Ax+Bv+Gw, \, \int_{0}^{T} Q^{-1}w(t).w(t)dt< \infty \}$. This space can be equipped with a quadratic norm based on $ \Pi_{0}^{-1}(x(0)).(x(0))+\int_{0}^{T} Q^{-1}w(t).w(t)dt+\int_{0}^{T} C^*R^{-1}Cx(t).x(t)dt$ and, once the affine term is removed, has a reproducing kernel Hilbert space (RKHS) structure. We refer to \cite{aubin2022Kalman} for more on this topic. Since we consider $\beta \neq \infty$, we authorize $x(.)$ to live beyond this RKHS. In other words the ``noise'' term $w(.)$ can be understood as a control, and by introducing $\beta$ we assume implicitly some extra noise on the model that was not present in \eqref{eq:2-1}. Moreover $\beta \rightarrow \infty$ implies that the minimizer $\hat w(.)$ is equal to $G^\ominus(\nicefrac{dx}{dt}-Ax(t)-Bv(t))$, with $G^\ominus$ the pseudo-inverse of $G$ for the Euclidean norm. Consequently the objective simplifies to
	\begin{align}
		\tilde J(A,B,x(.))&=\dfrac{\alpha}{2}\tr\left((A- A_{0})(A- A_{0})^{*}+(B- B_{0})(B- B_{0})^{*}\right)+\dfrac{1}{2} \Pi_{0}^{-1}(x(0)- x_{0}).(x(0)- x_{0})\nonumber \\
		&+\dfrac{1}{2}\int_{0}^{T} G^{\ominus,\top} Q^{-1}G^\ominus\left(\dfrac{dx}{dt}-Ax(t)-B v(t)\right).\left(\dfrac{dx}{dt}-Ax(t)-B v(t)\right)dt\nonumber \\
		&+\dfrac{1}{2}\int_{0}^{T} R^{-1}( y(t)- Cx(t)).( y(t)- Cx(t))dt \label{eq:2-6b}  ,
	\end{align}
	where it is effectively the noise $w(.)$ of \eqref{eq:2-1} that is penalized, and $	\tilde J$ corresponds to the traditional least square estimator used in Kalman smoothing \cite{Kailath2000,aubin2022Kalman}.
\end{remark}

\begin{remark}[Relation with EM] It is well-known that the Expectation-Maximization (EM) algorithm is an alternating minimization of a log-likelihood \cite{Neal1998}, which is the form of algorithm we propose in \Cref{sec:algo}. More precisely, given a probability space $(\cU,\tg)$, the relative entropy (Kullback--Leibler divergence) is defined as
	\begin{equation}\label{eq:kl_def}
		\KL(\mu|\tg)=\int_\cU \ln\left(\nicefrac{d\mu}{d\tg}(u)\right)d\mu(u)
	\end{equation}
	for $\mu$ absolutely continuous w.r.t.\ $\tg$ and $+\infty$ otherwise. In our case, we assume our observations $Y$ to be sampled according to $\tgY$ and $X$ serves as a latent, hidden random variable $X\in (\cX,\tg)$. We posit a joint distribution $p_{\qEM}(dx,dy)$ parametrized by an element $\qEM=(A,B)$ of the set $\Q=\cL(\R^{N};\R^{N})\times \cL(\R^{d};\R^{N})$. As presented in \cite{Neal1998}, the goal is to infer $\qEM$ by solving
	\begin{equation}
		\min_{\qEM\in\Q} \KL(\tgY|p_Y p_\qEM),
	\end{equation}
	where $p_Y p_{\qEM}(dy)=\int_{\cX} p_{\qEM}(dx,dy)$ is the marginal in $\cY$. The EM approach starts by minimizing a surrogate function of $\qEM$ upperbounding $\KL(\tgY|p_Y p_\qEM)$. For any $\pi\in\Pi(*,\tgY)=\{\pi\, | \, p_Y \pi = \tgY\}$, by the data processing inequality, i.e.\ KL of the marginals is smaller than KL of the plans,
	\[
	\KL(\tgY|p_Y p_\qEM)\le \KL(\pi|p_\qEM)=:L(\pi,\qEM).
	\]
	EM then proceeds by alternating minimizations of $L(\pi,\qEM)$ \cite[Theorem 1]{Neal1998}:
	\begin{align} 
		\label{eq:em-iterations-M-step}
		\qEM_{n} &= \argmin_{\qEM \in \Q} \KL(\pi_{n}|p_{\qEM}),\\
		\label{eq:em-iterations-E-step}
		\pi_{n+1}&=\argmin_{\pi\in\Pi(*,\tgY)}\KL(\pi|p_{\qEM_n}).
	\end{align}
	The above formulation consists in \eqref{eq:em-iterations-M-step}, optimizing  the parameters $\qEM_n$ at step $n$ (M-step), and then \eqref{eq:em-iterations-E-step}, optimizing  the joint distribution $\pi_{n+1}$ at step $n+1$ (E-step). This is actually what we propose as algorithm to minimize $J$ by minimizing alternatively in $X=(x(.),w(.))$ and $\theta=(A,B)$, $J$ being obtained as previously as the KL divergence of Gaussian measures. However making explicit the (Gaussian) measures underlying our parametrization goes beyond our scope and we now move to the study of our specific least-squares $J$.
\end{remark}

\section{Existence of a minimum and necessary condition
}\label{sec:existence}
Before searching for a minimum, we prove that $J$ has indeed one.
\begin{proposition}
	\label{prop3-1} The functional $J(Z)$ attains its infimum.
\end{proposition}

\begin{proof}
	The functional $J(Z)$ is continuous on $\cZ$ . It is also
	weakly lower semicontinuous. Indeed if $Z_{n}\rightharpoonup Z$ (weakly), then
	$A_{n}\rightarrow A$ in $\cL(\R^{N};\R^{N})$, $B_{n}\rightarrow B$ in $\cL(\R^{d};\R^{N})$, $x_{n}(.)\rightharpoonup x(.)$
	in $H^{1}(0,T;\R^{N})$, $w_{n}(.)\rightharpoonup w(.)$ in $L^{2}(0,T;\R^{m})$. We deduce that $(x_{n}(.))_n$ is equicontinuous, hence, by Ascoli's theorem, $x_{n}(.)\rightarrow x(.)$
	in $C^{0}([0,T];\R^{N})$.
	It follows that $A_{n}x_{n}(.)$$\rightarrow Ax(.)$ in $C^{0}([0,T];\R^{N})$
	and $\dfrac{dx_{n}}{dt}\rightharpoonup\dfrac{dx}{dt}$ in $L^{2}(0,T;\R^{N})$.
	Consequently we have 
	\begin{multline*}
		\dfrac{\alpha}{2}\tr\left((A_{n}-A_{0})(A_{n}-A_{0})^{*}+(B_{n}-B_{0})(B_{n}-B_{0})^{*}\right)\\+\dfrac{1}{2}\Pi_0^{-1}(x_{n}(0)-x_{0}).(x_{n}(0)-x_{0})+\dfrac{1}{2}\int_{0}^{T}R^{-1}(y(t)-Cx_{n}(t)).(y(t)-Cx_{n}(t))dt\\ \xrightarrow{n\rightarrow\infty} \dfrac{\alpha}{2}\tr\left((A-A_{0})(A-A_{0})^{*}+(B-B_{0})(B-B_{0})^{*}\right)\\+\dfrac{1}{2}\Pi_0^{-1}(x(0)-x_{0}).(x(0)-x_{0})+\dfrac{1}{2}\int_{0}^{T}R^{-1}(y(t)-Cx(t)).(y(t)-Cx(t))dt.
	\end{multline*}
	From the weak lower semicontinuity of the norm in the spaces $L^{2}(0,T;\R^{N})$
	and $L^{2}(0,T;\R^{m})$, we conclude easily that $J(Z)\leq\liminf\,J(Z_{n})$.
	If we consider a minimizing sequence $Z_{n}$, namely 
	\[
	J(Z_{n})\rightarrow\inf\,J(Z) \ge 0.
	\]
	Then, since $J(Z_{n})\leq J(0)$ for $n$ sufficiently large, it
	follows easily that the sequence $Z_{n}$ is bounded in $\cZ$.
	Since weakly closed bounded sets are weakly compact, we can extract a subsequence,
	still denoted $Z_{n}\rightharpoonup\widehat{Z}$ in $\cZ$
	weakly. From weak lower semicontinuity of $J$, we obtain $J$($\widehat{Z}$) $\leq\liminf$$\:J(Z_{n})=\inf\,J(Z)$.
	This implies that $\widehat{Z}$ is a minimum of $J(Z)$, which concludes
	the proof$.\blacksquare$ 
\end{proof}


We now check that $J(Z)$ has a Gâteaux differential in $Z$.
\begin{proposition}
	\label{prop3-2} The gradient $DJ(Z)\in\cZ$ is given by the formula 
	\begin{align}
		((DJ(Z),\widetilde{Z}))&=\tr\left(\left[\alpha(A-A_{0})+\int_{0}^{T}q(t)x^{*}(t)dt\right]\,\widetilde{A}^{*}\right)+\tr\left(\left[\alpha(B-B_{0})+\int_{0}^{T}q(t)v^{*}(t)dt\right]\,\widetilde{B}^{*}\right)\nonumber \\
		&+\left(\Pi_0^{-1}(x(0)-x_{0})+\int_{0}^{T}A^{*}q(t)-C^{*}R^{-1}(y(t)-Cx(t))dt\right).\widetilde{x}(0)\nonumber \\
		&+\int_{0}^{T}\left[-q(t)+\int_{t}^{T}A^{*}q(s)-C^{*}R^{-1}(y(s)-Cx(s))ds\right].\dfrac{d\widetilde{x}}{dt}(t)dt \nonumber \\
		&+\int_{0}^{T}(G^{*}q(t)+Q^{-1}w(t)).\widetilde{w}(t)dt\label{eq:3-1}
	\end{align}
	in which $Z=(A,B,x(.),w(.))$, $\widetilde{Z}=(\widetilde{A},\widetilde{B},\widetilde{x}(.),\widetilde{w}(.))$ and $q(t)$ is defined by 
	\begin{equation}
		q(t)=-\beta\left(\dfrac{dx}{dt}(t)-Ax(t)-Bv(t)-Gw(t)\right).\label{eq:3-2}
	\end{equation}
\end{proposition}

\begin{proof}
	From the definition of the Gâteaux differential in $Z$, we must check that  
	\begin{equation}
		\dfrac{d}{d\theta}J(Z+\theta\widetilde{Z})|_{\theta=0}=((DJ(Z),\widetilde{Z}))\label{eq:3-3}
	\end{equation}
	is equal to the right hand side of (\ref{eq:3-1}). We fix $Z$,
	$\widetilde{Z}$ with $q(t)$ defined by (\ref{eq:3-2}). As easily checked we can write 
	\begin{align}
		\dfrac{d}{d\theta}J(Z+\theta\widetilde{Z})|_{\theta=0}&=\tr\left(\left[\alpha(A-A_{0})+\int_{0}^{T}q(t)x^{*}(t)dt\right]\,\widetilde{A}^{*}\right)\nonumber \\
		&+\tr\left(\left[\alpha(B-B_{0})+\int_{0}^{T}q(t)v^{*}(t)dt\right]\,\widetilde{B}^{*}\right)\nonumber \\
		&+\Pi_0^{-1}(x(0)-x_{0}).\widetilde{x}(0)-\int_{0}^{T}q(t).\left(\dfrac{d}{dt}\widetilde{x}-A\widetilde{x}(t)\right)dt\nonumber \\
		&-\int_{0}^{T}C^{*}R^{-1}(y(t)-Cx(t)).\widetilde{x}(t)dt+\int_{0}^{T}(G^{*}q(t)+Q^{-1}w(t)).\widetilde{w}(t)dt. \label{eq:3-4}
	\end{align}
	We then replace $\widetilde{x}(t)$ with $\widetilde{x}(0)+\int_{0}^{t}\frac{d}{ds}\widetilde{x}(s)ds$.
	We perform a change of integration and some rearrangements to obtain the relation (\ref{eq:3-1}).
\end{proof}
If $\widehat{Z}=(\widehat{A},\widehat{B},\widehat{x}(.),\widehat{w}(.))$ is a point of minimum for $J(Z)$, it follows from formula (\ref{eq:3-1}) that the corresponding $\widehat{q}(t)$ defined by 
\begin{equation}
	\widehat{q}(t)=-\beta\left(\dfrac{d\widehat{x}}{dt}(t)-\widehat{A}\widehat{x}(t)-\widehat{B}v(t)-G\widehat{w}(t)\right)\label{eq:3-5}
\end{equation}
satisfies 
\begin{align*}
	\Pi_0^{-1}(\widehat{x}(0)-x_{0})+\int_{0}^{T}(A^{*}\widehat{q}(t)-C^{*}R^{-1}(y(t)-C\widehat{x}(t))dt&=0\\
	-\widehat{q}(t)+\int_{t}^{T}(A^{*}\widehat{q}(s)-C^{*}R^{-1}(y(s)-C\widehat{x}(s))ds&=0,\forall t.
\end{align*}
It follows that $\dfrac{d\widehat{q}}{dt}$ is well defined. Differentiating the previous equation, and reordering \eqref{eq:3-5}, we obtain the following system of optimality conditions for \eqref{eq:3-1}
\begin{align}
	\dfrac{d\widehat{x}}{dt}&=\widehat{A}\widehat{x}(t)+\widehat{B}v(t)-(GQ^{-1}G^{*}+\dfrac{I}{\beta})\widehat{q}(t), \, &&\widehat{x}(0)=x_{0}-\Pi_0\widehat{q}(0), \label{eq:3-6}\\
	-\dfrac{d\widehat{q}}{dt}&=\widehat{A}^{*}\widehat{q}(t)-C^{*}R^{-1}C(y(t)-C\widehat{x}(t)),\, &&\widehat{q}(T)=0, \nonumber\\
	&\alpha(\widehat{A}-A_{0})+\int_{0}^{T}\widehat{q}(t)\widehat{x}^{*}(t)dt=0,\label{eq:3-7}\\
	&\alpha(\widehat{B}-A_{0})+\int_{0}^{T}\widehat{q}(t)v^{*}(t)dt=0, \nonumber
\end{align}
with $\widehat{w}(t)$ given by 
\begin{equation}
	\widehat{w}(t)=-QG^{*}\widehat{q}(t).\label{eq:3-8}
\end{equation}
Note that \eqref{eq:3-7} has an interesting structure, decomposing the optimal $\widehat{A}$ (resp.\ $\widehat{B}$) as a sum of rank 1 tensor products between the covector $q$ and the trajectory $x$ (resp.\ covector $q$ and control $v$).
\section{Alternating minimization algorithm}\label{sec:algo}

The relations (\ref{eq:3-6}), (\ref{eq:3-7}) can be interpreted
as a fixed point problem for the pair $(\widehat{A},\widehat{B}$).
If we fix the pair $(\widehat{A},\widehat{B})$ then we obtain the
pair $(\widehat{x}(.),\widehat{q}(.))$ by solving the system of forward
backward equations (\ref{eq:3-6}). Next for fixed $(\widehat{x}(.),\widehat{q}(.))$
we obtain $(\widehat{A},\widehat{B})$ by the formulas (\ref{eq:3-7}). This corresponds also to an alternating minimization of $J$, which happens in two steps.

The first part is associated to a control problem, formulated
as a calculus of variations problem
\begin{align}
	\min_{x(.),w(.)}K(\widehat{A},\widehat{B};x(.),w(.))&:=\dfrac{1}{2}\Pi_0^{-1}(x(0)-x_{0}).(x(0)-x_{0})\nonumber\\
	&+\dfrac{\beta}{2}\int_{0}^{T}\left|\dfrac{dx}{dt}-\widehat{A}x(t)-\widehat{B}v(t)-Gw(t)\right|^{2}dt \nonumber\\ &+\dfrac{1}{2}\int_{0}^{T}Q^{-1}w(t).w(t)dt+\dfrac{1}{2}\int_{0}^{T}R^{-1}(y(t)-Cx(t)).(y(t)-Cx(t))dt.\label{eq:3-9} 
\end{align}
The second part is associated to an optimization problem 
\begin{align}
	\min_{A,B}L(A,B;\widehat{x}(.),\widehat{w}(.))&:=\dfrac{\alpha}{2}\tr\left((A-A_{0})(A-A_{0})^{*}+(B-B_{0})(B-B_{0})^{*}\right)\nonumber \\&+\dfrac{\beta}{2}\int_{0}^{T}\left|\dfrac{d\widehat{x}}{dt}-A\widehat{x}(t)-Bv(t)-G\widehat{w}(t)\right|^{2}dt.\label{eq:3-11}
\end{align}
It is important to notice that the two problems (\ref{eq:3-9}), (\ref{eq:3-11}) are convex quadratic and have a unique solution, whereas the original problem (\ref{eq:2-6}) is not convex. This highlights the usefulness of algorithm that we propose to find a local optimum $\widehat{Z}$ of $J(Z)$.


We initialize the algorithm with $(A_{0},B_{0})$. For $n\ge 0$, knowing $A_{n},B_{n}$ we
define uniquely the pair $(x_{n}(.),w_{n}(.))$ which minimizes $K(A_{n},B_{n};x(.),w(.))$.
This leads immediately to the existence and uniqueness of the pair
$x_{n}(.),q_{n}(.)$ solution of the system of forward-backward relations 
\begin{align}
	\dfrac{dx_{n}}{dt}&=A_{n}x_{n}(t)+B_{n}v(t)-\left(GQ^{-1}G^{*}+\dfrac{I}{\beta}\right)q_{n}(t), \, &&x_{n}(0)=x_{0}-\Pi_0q_{n}(0)\label{eq:4-1}\\
	-\dfrac{dq_{n}}{dt}&=(A_{n})^{*}q_{n}(t)-C^{*}R^{-1}C(y(t)-Cx_{n}(t)), \, &&q_{n}(T)=0,
\end{align}
with $w_{n}(t)$ given by 
\begin{equation}
	w_{n}(t))=-QG^{*}q_{n}(t). \label{eq:4-2}
\end{equation}
We then define $A_{n+1},B_{n+1}$ by minimizing $L(A,B;x_{n}(.),w_{n}(.))$.
We obtain 
\begin{align}
	&\alpha(A_{n+1}-A_{0})-\beta\int_{0}^{T}\left(\dfrac{dx_{n}}{dt}-A_{n+1}x_{n}(t)-B_{n+1}v(t)-Gw_{n}(t)\right)(x_{n}(t))^{*}dt=0, \label{eq:4-3}\\
	&\alpha(B_{n+1}-B_{0})-\beta\int_{0}^{T}\left(\dfrac{dx_{n}}{dt}-A_{n+1}x_{n}(t)-B_{n+1}v(t)-Gw_{n}(t)\right)(v(t))^{*}dt=0.\nonumber
\end{align}
Using (\ref{eq:4-1})-\eqref{eq:4-2}, the term in $w_n$ canceling out with one of those in $q_n$, we can rewrite the equation (\ref{eq:4-3}) as follows by factorizing
\begin{multline}
	A_{n+1}\left(\alpha I+\beta\int_{0}^{T}x_{n}(t)(x_{n}(t))^{*}dt\right)+\beta B_{n+1}\int_{0}^{T}v(t)(x_{n}(t))^{*}dt\\=\alpha A_{0}+\beta A_{n}\int_{0}^{T}x_{n}(t)(x_{n}(t))^{*}dt+\beta B_{n}\int_{0}^{T}v(t)(x_{n}(t))^{*}dt-\int_{0}^{T}q_{n}(t)(x_{n}(t))^{*}dt,\label{eq:4-4}
\end{multline}
\begin{multline*}
	\beta A_{n+1}\int_{0}^{T}x_{n}(t)(v(t))^{*}dt+B_{n+1}(\alpha I+\beta\int_{0}^{T}v(t)(v(t))^{*}dt)\\=\alpha B_{0}+\beta A_{n}\int_{0}^{T}x_{n}(t)(v(t))^{*}dt+\beta B_{n}\int_{0}^{T}v(t)(v(t))^{*}dt-\int_{0}^{T}q_{n}(t)(v(t))^{*}dt.
\end{multline*}

Our main result is the convergence of the alternating minimization scheme to extremal points for the first-order optimality conditions.
\begin{theorem}
	\label{theo4-1} The sequence $J(Z^{n})$ is decreasing. The sequence
	$Z^{n}$ is bounded and $Z^{n+1}-Z^{n}\rightarrow0$. Limits of converging
	subsequences of $Z^{n}$ are solutions of the set of necessary conditions
	(\ref{eq:3-6}), (\ref{eq:3-7}). 
\end{theorem}
\begin{proof}
	We compute the two differences $K(A_{n+1},B_{n+1};x_{n}(.),w_{n}(.))-K(A_{n+1},B_{n+1};x_{n+1}(.),w_{n+1}(.))$
	and $L(A_{n},B_{n};x_{n}(.),w_{n}(.))-L(A_{n+1},B_{n+1};x_{n}(.),w_{n}(.))$
	which are nonnegative numbers, since $(x_{n+1}(.),w_{n+1}(.))$ minimizes
	$K(A_{n+1},B_{n+1};x(.),w(.))$ and $(A_{n+1},B_{n+1})$ minimizes $L(A,B;x_{n}(.),w_{n}(.))$. Since we know the minimizers, we can use the fact that a quadratic function $f(z)$, with Hessian $H$ and minimum $\bar z$, satisfies $f(z)-f(\bar z)=\frac{1}{2}(z-\bar z)^* H (z-\bar z)$. This is a completion of square argument. Consequently, we have 
	\begin{align}
		&K(A_{n+1},B_{n+1};x_{n}(.),w_{n}(.))-K(A_{n+1},B_{n+1};x_{n+1}(.),w_{n+1}(.))\nonumber \\&=\dfrac{1}{2}\Pi_0^{-1}(x_{n}(0)-x_{n+1}(0)).(x_{n}(0)-x_{n+1}(0))\nonumber \\
		&+\dfrac{\beta}{2}\int_{0}^{T}\left|\dfrac{d}{dt}(x_{n}-x_{n+1})-A_{n+1}(x_{n}(t)-x_{n+1}(t))-G(w_{n}(t)-w_{n+1}(t))\right|^{2}dt\nonumber \\ &+\dfrac{1}{2}\int_{0}^{T}Q^{-1}(w_{n}(t)-w_{n+1}(t)).(w_{n}(t)-w_{n+1}(t))dt\nonumber \\ &+\dfrac{1}{2}\int_{0}^{T}R^{-1}C(x_{n}(t)-x_{n+1}(t)).C(x_{n}(t)-x_{n+1}(t))dt.\label{eq:4-5}
	\end{align}
	Similarly we can write also 
	\begin{align}
		&L(A_{n},B_{n};x_{n}(.),w_{n}(.))-L(A_{n+1},B_{n+1};x_{n}(.),w_{n}(.))\nonumber \\&=\dfrac{\alpha}{2}\tr\:\left((A_{n}-A_{n+1})(A_{n}-A_{n+1})^{*}+(B_{n}-B_{n+1})(B_{n}-B_{n+1})^{*}\right)\nonumber \\&+\dfrac{\beta}{2}\int_{0}^{T}\left|(A_{n}-A_{n+1})x_{n}(t)+(B_{n}-B_{n+1})v(t)\right|^{2}dt.\label{eq:4-6}
	\end{align}
	Adding up, we see that the left hand side is $J(Z^{n})-J(Z^{n+1})$, so we obtain 
	\begin{align}
		J(Z^{n})-J(Z^{n+1})
		&=\dfrac{1}{2}\Pi_0^{-1}(x_{n}(0)-x_{n+1}(0)).(x_{n}(0)-x_{n+1}(0))\nonumber \\&+\dfrac{\beta}{2}\int_{0}^{T}\left|\dfrac{d}{dt}(x_{n}-x_{n+1})-A_{n+1}(x_{n}(t)-x_{n+1}(t))-G(w_{n}(t)-w_{n+1}(t))\right|^{2}dt\nonumber \\
		&+\dfrac{1}{2}\int_{0}^{T}Q^{-1}(w_{n}(t)-w_{n+1}(t)).(w_{n}(t)-w_{n+1}(t))dt\nonumber \\
		&+\dfrac{1}{2}\int_{0}^{T}R^{-1}C(x_{n}(t)-x_{n+1}(t)).C(x_{n}(t)-x_{n+1}(t))dt\nonumber \\ &+\dfrac{\alpha}{2}\tr\:\left((A_{n}-A_{n+1})(A_{n}-A_{n+1})^{*}+(B_{n}-B_{n+1})(B_{n}-B_{n+1})^{*}\right)\nonumber  \\&+\dfrac{\beta}{2}\int_{0}^{T}\left|(A_{n}-A_{n+1})x_{n}(t)+(B_{n}-B_{n+1})v(t)\right|^{2}dt, \label{eq:4-7}
	\end{align}
	which is a nonnegative quantity. It follows that the sequence $J(Z^{n})$ is decreasing. Since it
	is nonnegative, it converges. From the relation (\ref{eq:4-7})
	we see that $Z^{n}-Z^{n+1}\rightarrow0$ in $\|.\|_\cZ$. Since $J(Z^{n})\leq J(Z^{0})$,
	the sequence $Z^{n}$ is bounded in $\cZ$. If we extract a subsequence  which converges weakly to $\widehat{Z}$, also
	noted $Z^{n}$ without loss of generality, then $A_{n}\rightarrow\widehat{A},B_{n}\rightarrow\widehat{B}$ and
	$x_{n}(.)\rightarrow\widehat{x}(.)$ in $H^{1}(0,T;\R^{N})$ weakly,
	hence strongly in $C^{0}([0,T];\R^{N})$. From (\ref{eq:4-1}), $q_{n}(.)\rightarrow\widehat{q}(,)$
	in $H^{1}(0,T;\R^{N})$ weakly and strongly in $C^{0}([0,T];\R^{N})$.
	Therefore, from (\ref{eq:4-2}) $w_{n}(.)\rightarrow\widehat{w}(.)$
	in $L^{2}(0,T;\R^{N})$. Finally $Z_{n}\rightarrow\widehat{Z.}$ We
	can thus take the limit in equations (\ref{eq:4-1}), (\ref{eq:4-4})
	and obtain that $\widehat{Z}$ is solution of the set of equations
	(\ref{eq:3-6}), (\ref{eq:3-7}). This concludes the proof.
\end{proof}

\vskip 6mm
\noindent{\bf Acknowledgements}

\noindent
The first author was funded by the European Research Council (grant REAL 947908). The second author was supported by the National Science Foundation under grants NSF-DMS-1905449, NSF-DMS-2204795 and grant from the SAR Hong Kong RGC GRF 14301321.      The third author acknowledges partial financial support for this work from  a General Research Fund by the Research Grants Council (RGC) of Hong Kong SAR, China (Project No. 11303421), a grant from ITF - Guangdong-Hong Kong Technology Cooperation Funding Scheme (Project Ref. No. GHP/145/20), and a Math and Application Project (2021YFA1003504) under the National Key R\&D Program.


\begin{thebibliography}{99}

	\providecommand{\natexlab}[1]{#1}
	\providecommand{\url}[1]{\texttt{#1}}
	\expandafter\ifx\csname urlstyle\endcsname\relax
	\providecommand{\doi}[1]{doi: #1}\else
	\providecommand{\doi}{doi: \begingroup \urlstyle{rm}\Url}\fi

	\bibitem[1a]{leger2023gradient}
	Pierre-Cyril Aubin-Frankowski and Flavien L{\'e}ger.
	\newblock Gradient descent with a general cost.
	\newblock 2023
	\newblock \url{https://arxiv.org/abs/2305.04917}.

	
	\bibitem[1b]{aubin2022Kalman}
	Pierre-Cyril Aubin-Frankowski and Alain Bensoussan.
	\newblock The reproducing kernel {Hilbert} spaces underlying linear {SDE}
	estimation, {Kalman} filtering and their relation to optimal control.
	\newblock \emph{Pure and Applied Functional Analysis}, 2023.
	\newblock (to appear) \url{https://arxiv.org/abs/2208.07030}.
			
	\bibitem[2]{Bensoussan2018}
	Alain Bensoussan.
	\newblock \emph{Estimation and Control of Dynamical Systems}.
	\newblock Springer International Publishing, 2018.
	
	\bibitem[3]{Bensoussan2021}
	Alain Bensoussan, Fatih Gelir, Viswanath Ramakrishna, and Minh-Binh Tran.
	\newblock Identification of linear dynamical systems and machine learning.
	\newblock \emph{Journal of Convex Analysis}, 28\penalty0 (2):\penalty0
	311--328, 2021.
	
	\bibitem[4]{DaPrato2006}
	Giuseppe Da~Prato.
	\newblock \emph{An Introduction to Infinite-Dimensional Analysis}.
	\newblock Springer Berlin Heidelberg, 2006.
	\newblock \doi{10.1007/3-540-29021-4}.
	
	
	\bibitem[5]{Dashti2013}
	Masoumeh Dashti, Kody J~H Law, Andrew~M Stuart, and Jochen Voss.
	\newblock Map estimators and their consistency in bayesian nonparametric
	inverse problems.
	\newblock \emph{Inverse Problems}, 29\penalty0 (9):\penalty0 095017, Sep 2013.
	\newblock ISSN 0266-5611, 1361-6420.
	\newblock \doi{10.1088/0266-5611/29/9/095017}.
		
	\bibitem[6]{Guth2022}
	Philipp~A. Guth, Claudia Schillings, and Simon Weissmann.
	\newblock Ensemble Kalman filter for neural network based one-shot inversion.
	\newblock In \emph{Optimization and Control for Partial Differential
		Equations}, pages 393--418. De Gruyter, 2022.
	\newblock \doi{10.1515/9783110695984-014}.
	
	\bibitem[7]{Kailath2000}
	Thomas Kailath, Ali~H Sayed, and Babak Hassibi.
	\newblock \emph{Linear Estimation}.
	\newblock Prentice Hall information and system sciences series. Pearson, 2000.
	
	\bibitem[8]{Neal1998}
	Radford~M. Neal and Geoffrey~E. Hinton.
	\newblock A view of the {EM} algorithm that justifies incremental, sparse, and
	other variants.
	\newblock In \emph{Learning in Graphical Models}, pages 355--368. Springer
	Netherlands, 1998.
	\newblock \doi{10.1007/978-94-011-5014-9_12}.
		
	\bibitem[9]{Qin2022}
	S.~Joe Qin.
	\newblock Latent vector autoregressive modeling and feature analysis of high
	dimensional and noisy data from dynamic systems.
	\newblock \emph{{AIChE} Journal}, 68\penalty0 (6), April 2022.
	\newblock \doi{10.1002/aic.17703}.
	
	\bibitem[10]{Qin2020}
	S.~Joe Qin, Yining Dong, Qinqin Zhu, Jin Wang, and Qiang Liu.
	\newblock Bridging systems theory and data science: A unifying review of
	dynamic latent variable analytics and process monitoring.
	\newblock \emph{Annual Reviews in Control}, 50:\penalty0 29--48, 2020.
	\newblock \doi{10.1016/j.arcontrol.2020.09.004}.
	
	\bibitem[11]{Yu2022}
	Jiaxin Yu and S.~Joe Qin.
	\newblock Latent state space modeling of high-dimensional time series with a
	canonical correlation objective.
	\newblock \emph{{IEEE} Control Systems Letters}, 6:\penalty0 3469--3474, 2022.
	\newblock \doi{10.1109/lcsys.2022.3183895}.




\end{thebibliography}
\end{document}